\documentclass{amsart}

\usepackage{amssymb}
\usepackage{amsmath}
\usepackage{amsthm}

\usepackage{graphicx}

\newtheorem {theorem}{Theorem}
\newtheorem {proposition}[theorem]{Proposition}

\newtheorem {lemma}[theorem]{Lemma}
\theoremstyle{definition}
\newtheorem{df}{Definition}

\theoremstyle{remark}

\numberwithin{equation}{section}

\newcommand{\al}{\alpha}

\newcommand{\er}{\varepsilon}
\newcommand{\za}{\zeta}

\newcommand{\we}{w}

\newcommand{\ph}{\varphi}

\newcommand{\Rbb}{\mathbb R}

\newcommand{\qle}{\ell}
\newcommand{\Om}{\Omega}
\newcommand{\om}{\omega}
\newcommand{\omt}{\widetilde{\omega}}
\newcommand{\ccom}{\mathcal{C}^{1,\omega}}
\newcommand{\Cbb}{\mathbb C}
\newcommand{\bmo}{\mathrm{BMO}}
\newcommand{\mo}{\mathrm{MO}}

\newcommand{\har}{\chi}
\newcommand{\beu}{B_\Omega}
\newcommand{\bloch}{\mathcal B}
\newcommand{\fti}{\widetilde{f}}
\newcommand{\qti}{\widetilde{Q}}
\newcommand{\Tbb}{\mathbb T}
\newcommand{\Dbb}{\mathbb D}
\newcommand{\qom}{Q|\Om}
\newcommand{\ext}{\varkappa}

\begin{document}

\author[E.~Doubtsov]{Evgueni Doubtsov}

\address{St.~Petersburg Department
of V.A.~Steklov Mathematical Institute,
Fontanka 27, St.~Petersburg
191023, Russia}

\email{dubtsov@pdmi.ras.ru}

\title
{Restricted Beurling transforms on Campanato spaces}

\author[A.~V.~Vasin]{Andrei V.~Vasin}

\address{State University of Maritime and Inland Shipping,
Dvinskaya st.~5/7, St.~Petersburg 198035, Russia}

\email{andrejvasin@gmail.com}

\thanks{The authors were supported by RFBR (grant No.~14-01-00198a).}


\date{}


\begin{abstract}
Let $\Om\subset\Cbb$ be a bounded domain with $\ccom$-smooth boundary, where
$\om$ is a Dini-smooth modulus of continuity.
We prove that the restricted Beurling transform is bounded on the Campanato space $\bmo_\om(\Om)$.
\end{abstract}

\maketitle

\section{Introduction}\label{s_int}
\textsl{The Beurling transform}
is the principal value convolution operator
\[
{B}f(z) = -\frac{1}{\pi} \textrm{p.v.} \int_{\Cbb} f(z-w)\frac{1}{w^2}\, dw.
\]
Given a bounded domain $\Om\subset\Cbb$,
in the present paper we consider the corresponding modification of ${B}$.
Namely, \textsl{the restricted Beurling transform} $\beu$ is defined as
\[
\beu f(z) = {B}(\chi_\Om f)(z), \quad z\in \Om,
\]
where $\chi_\Om$ is the characteristic function of $\Om$.

\subsection{Motivation: applications to quasiregular mappings}\label{ss_classlip}
The studies of $\beu$ are motivated by applications
in the theory of quasiregular mappings (see \cite{AIM09}).
In the setting of the classical Lipschitz spaces
$\Lambda^\al(\Om)$, $0<\al< 1$,
the following result is known to be crucial; see \cite{MOV09}.

\begin{theorem}[see {\cite[Main Lemma]{MOV09}}]\label{t_mov}
Let $\Om\subset\Cbb$ be a bounded domain with $\mathcal{C}^{1+\al}$-smooth boundary, $0<\al< 1$.
Then the restricted Beurling transform $\beu$ is bounded on $\Lambda^\al(\Om)$.
\end{theorem}

In particular, the above theorem is applied in \cite{MOV09} to study the principal solutions
of the Beltrami equation
\[
\overline{\partial}f = \mu\partial f,\quad \textrm{a.e.~on\ } \Om,
\]
under assumption that the Beltrami coefficient $\mu$, $\|\mu\|_{L^\infty(\Om)} < 1$,
is in $\Lambda^\al(\Om)$, $0<\al< 1$.
See also \cite{CMO13}, where the Lipschitz conditions are replaced by Sobolev or Besov conditions.

Motivated by the applications mentioned above, we are looking for extensions of Theorem~\ref{t_mov}
to larger classes of bounded domains $\Om\subset\Cbb$.
A more precise problem is formulated as follows:
Under weaker restrictions on $\partial\Om$, find $\beu$-invariant spaces $X(\Om)$, $X(\Om)\subset L^\infty(\Om)$.
To solve the problem, we consider spaces defined by means of regular moduli of continuity.

\subsection{Moduli of continuity}
\begin{df}
  An increasing continuous
  function $\om: [0, +\infty)\to [0, +\infty)$, $\om(0)=0$,
  is called \textsl{a modulus of continuity}.
\end{df}

The Lipschitz space $\Lambda^{\om}(\Om)$ consists
  of those functions $f:\Om\to\Cbb$ for which there exists a constant $C>0$ such that
  \[
  |f(z)-f(w)| \le C \om(|z-w|),\quad z,w\in \Om.
  \]
For an interval $(a,b)\subset\Rbb$, the space $\Lambda^\om(a, b)$ is defined analogously.
Also, we write $\ph\in\mathcal{C}^{1,\om}(a, b)$ if $\ph$ is differentiable on $(a, b)$
 and $\ph^\prime\in \Lambda^\om(a, b)$.

A modulus of continuity $\om$ is called \textsl{Dini-smooth} if the integral
\[
\int_0 \frac{\om(t)}{t}\, dt
\]
converges.

Also, we use the following regularity condition: there exists $\er\in (0,1)$ such that
the quotient $\om(t)/{t^\er}$ is \textsl{almost decreasing}, that is,
\begin{equation}\label{e_almdec}
\frac{\om(t)}{t^\er} \le C \frac{\om(s)}{s^\er}, \quad t>s>0,
\end{equation}
for a universal constant $C>0$.

\begin{df}
A modulus of continuity $\om$ is called \textsl{regular} if $\om$ is Dini-smooth and
property \eqref{e_almdec} holds.
\end{df}
In what follows, we assume that $\om$ is regular if not otherwise stated.
For $\beta>0$, the logarithmic function
\[
\left(\log\frac{e}{t} \right)^{-1-\beta}, \quad 0<t<1,
\]
may serve as a working example of a regular modulus of continuity restricted to $(0,1)$.

\subsection{Restricted Beurling transform on $\Lambda^\om(\Om)$}
Given a modulus of continuity $\om$, we say that $\Om\subset\Cbb$ is a domain
with $\ccom$-smooth boundary if $\partial\Om$ is a $\mathcal{C}^1$ curve
whose unit normal vector is in $\Lambda^\om$ as a function on the curve.
An equivalent technical reformulation of this assumption is given in Section~\ref{ss_r0}.

\begin{proposition}[{\cite[Theorem~1]{Va15}}]\label{p_lip}
Let $\om$ be a Dini-smooth
modulus of continuity and let $\Om\subset\Cbb$ be a bounded domain with $\ccom$-smooth boundary.
Then $\beu$ maps $\Lambda^{\om}(\Om)$ into $\Lambda^{\omt}(\Om)$, where
\[
\omt(x) = \int_0^x \frac{\om(t)}{t}\, dt + x \int_x^1 \frac{\om(t)}{t^2}\, dt
\]
is the conjugate modulus of continuity.
\end{proposition}

As shown in \cite{Va15}, Proposition~\ref{p_lip} is, in a sense, sharp in the scale of Lipschitz spaces.
Moreover, general {C}alder\'on--{Z}ygmund operators
map $\Lambda^{\om}$ into $\Lambda^{\omt}$ even in the case, where the smoothness of $\partial\Om$
does not affect the computations (see, for example, \cite{KK13} for various results of this type).
Also, we have $\omt_\al=\om_\al$ for $\om_\al(t)=t^\al$, $0<\al< 1$.
So, on the one hand, Proposition~\ref{p_lip} is a sharp and direct extension of Theorem~\ref{t_mov}.
On the other hand, Proposition~\ref{p_lip} does not solve the $\beu$-invariance problem formulated in Section~\ref{ss_classlip}.
In the present paper, we show that an appropriate choice for a $\beu$-invariant space $X(\Omega)\subset L^\infty(\Om)$
is the Campanato space $\bmo_\om(\Om)$ with a regular $\om$.

\subsection{Restricted Beurling transform on Campanato spaces}
In what follows, $Q\subset\Cbb$
denotes a square
with edges parallel to the coordinates axes, $\qle$ denotes the side length of $Q$,
and $|Q|=\qle^2$ is the area of $Q$.
Let $dA$ denote the area measure on $\Cbb$. Given $1\le p < \infty$ and a modulus of continuity $\om$,
the Campanato space $\bmo_{\om} = \bmo_{\om, p} (\Cbb)$ consists of those $g\in L^p_{loc}(\Cbb)$ for which
\begin{equation}\label{e_bmo_df}
\|g\|_{\bmo_{\om, p}} = \sup_{Q\subset \Cbb} \frac{1}{\om(\qle)} \|g - g_Q\|_{L^p(Q, dA/|Q|)}
< \infty,
\end{equation}
where
\[
g_Q = \frac{1}{|Q|} \int_Q g(z)\, dA(z)
\]
is the standard integral mean of $g$ over $Q$.
In fact, the arguments used in the studies of the classical space $\bmo(\Rbb^n)$, $n\ge 1$,
guarantee that, for all $1\le p < \infty$, the seminorms under consideration define the same space
(see, for example, \cite{MMNO2k}, where methods from \cite{CRT08} are adapted).
We use the equivalence of the corresponding seminorms in the proof of Theorem~\ref{t_main};
see Section~\ref{ss_sterm}.
So, below we consider the seminorm with $p=1$ and we write $\bmo_{\om}$ in the place of $\bmo_{\om, 1}$.
The space $\bmo_\om$ is a generalization of $\Lambda^\al$, $0<\al<1$:
for $\om_\al(t) = t^\al$, one has $\bmo_{\om_\al} = \Lambda^\al$ (see \cite{Ca63, Me64}).
In fact,
$\bmo_\om \subset \Lambda^{\omt}$ for all Dini-regular $\om$; see \cite{Sp65}.
Also, the spaces $\bmo_{\om_\al} = \Lambda^\al$, $0<\al< 1$, are known to be invariant for certain Calder\'on--Zygmund
convolution operators
(see \cite{Pe66}).

To define $\bmo_\om(\Om)$ for a domain $\Om\subset \Cbb$, we use formula~\eqref{e_bmo_df},
replacing the norm in $L^p(Q, dA/|Q|)$ by that in $L^p(Q\cap \Om, dA/|Q|)$ with $p=1$.
So, $\bmo_\om(\Om)$ consists of those $f\in L^\infty(\Om)$ for which
\begin{equation}\label{e_camp_df}
\|f\|_{\bmo_{\om}(\Om)} = \sup_{Q\subset \Cbb} \frac{1}{\om(\qle)}
\frac{1}{|Q|} \int_{Q\cap\Om} |f(z) - f_{Q|\Om}|\, dA(z)
< \infty,
\end{equation}
where
\[
f_{Q|\Om} = \frac{1}{|Q\cap\Om|} \int_{Q\cap\Om} f(z)\, dA(z)
\]
is the integral mean of $f$ over $Q\cap\Om$.

Observe that in the case of the classical space $\bmo=\bmo(\Cbb)$, the above
definition reduces to that of $\bmo^r(\Om)$, \textsl{the restricted $\bmo$} on $\Om$
(see, for example, \cite{Jo80}).
While several $\bmo$ space are known for bounded domains, $\bmo^r(\Om)$
is usually considered as a standard one.
For regular $\om$, we have $\bmo_\om \subset \Lambda^{\omt}$, hence,
on the one hand,
$g\in\bmo_\om$ implies that the restriction $g\left|_{\Om}\right.$
is in $L^\infty(\Om)$ and $\|g\left|_{\Om}\right.\|_{\bmo_\om(\Om)} < \infty$.
On the other hand, if $f\in\bmo_\om(\Om)$, then Lemma~\ref{l_ext} provides an
extension $\fti\in\bmo_\om$ such that $\fti\left|_{\Om}\right. = f$.
Hence, $\bmo_\we(\Om)$ is a weighted restricted $\bmo$ on $\Om$.
Also, the property $\bmo_\om(\Om)\subset L^\infty(\Om)$ is a natural assumption
in the definition of $\bmo_\om(\Om)$.

The main result of the present paper guarantees that the space $\bmo_\om(\Omega)$
is $\beu$-invariant for all domains $\Om$ and moduli of continuity $\om$ under consideration.

\begin{theorem}\label{t_main}
Let $\om$ be a regular modulus of continuity and let
$\Om\subset\Cbb$ be a bounded domain with $\ccom$-smooth boundary.
Then the restricted Beurling transform $\beu$ is bounded on $\bmo_\om(\Om)$.
\end{theorem}

\subsection*{Notation}
The symbol $C$, with or without subscripts, denotes an absolute constant
whose value may change from line to line.
We write $E \lesssim F$ if $E \le CF$ for a constant $C>0$;
also, $E \lesssim F \lesssim E$ is abreviated as $E \approx F$.

\textsl{A square} on the complex plane $\Cbb$ is always a square with edges parallel to the coordinate axes.
Given a square $Q$ and a constant $a>0$, $aQ$ denotes the square with side length $a\qle(Q)$ and concentric with $Q$.

\subsection*{Organization of the paper}
Auxiliary facts are collected in Section~\ref{s_aux};
in particular, we prove embedding and extension results
(Lemmas~\ref{l_emb} and \ref{l_ext}, respectively).
The main result, Theorem~\ref{t_main}, is proved in Section~\ref{s_main}.

\section{Auxiliary results}\label{s_aux}

\subsection{Regular moduli of continuity}
Let $\om$ be a regular modulus of continuity.
First, property~\eqref{e_almdec} trivially implies that
the quotient $\om(t)/t$ is almost decreasing, that is,
there exists a constant $C>0$ such that
\[
\frac{\om(t)}{t} \le C \frac{\om(s)}{s} \quad \textrm{for all\ } t>s>0.
\]
Secondly, \eqref{e_almdec} is known to be equivalent the following integral
condition:
\begin{equation}\label{e_weak}
x\int_x^1 \frac{\om(t)}{t^2}\, dt \le C\om(x)
\end{equation}
for a positive constant $C$; see, for example, \cite{Ja76, Ja80}.

\subsection{$\ccom$-smooth boundary: technical properties and assumptions}\label{ss_r0}
Standard geometric arguments guarantee (see, for example, \cite{MOV09, To13} for related
conditions in $\Rbb^n$, $n\ge 2$) that $\Om\subset\Cbb$ is a bounded domain with a $\ccom$-smooth boundary
if and only if for each point $\za\in\partial\Om$, there exists
a square $Q=Q(\za)$ such that $\Om\cap Q(\za)$ is a subgraph of a $\mathcal{C}^{1,\om}$ function.
More precisely, given $\za\in\partial\Om$, we may assume
(probably, after a suitable rotation by $\pm\frac{\pi}{2}$ or $\pi$)
that the slope of the tangent to $\partial\Om$ at $\za$ is positive.
If the slope of the tangent is at least $1/10$, then consider a sufficiently small square $Q=Q(\za)$
such that $\za$ is the lower left vertex of $Q$.
If $\partial\Om$ intersects the right edge of $Q$, then let $I(\za)$ denote the lower edge of $Q$;
if $\partial\Om$ intersects the upper edge of $Q$, then let $I(\za)$ denote the right edge of $Q$.
Below we consider the first case; in the second case, we argue analogously after a suitable rotation.
Also, if the slope of the tangent at $\za$ is less than $1/10$, then it suffices to consider
a sufficiently small square $Q$ such that $\za$ is the center of the left edge of $Q$ and to define $I(\za)$ as the lower
edge of $Q$.

\begin{figure}
  \includegraphics[width=0.5\textwidth]{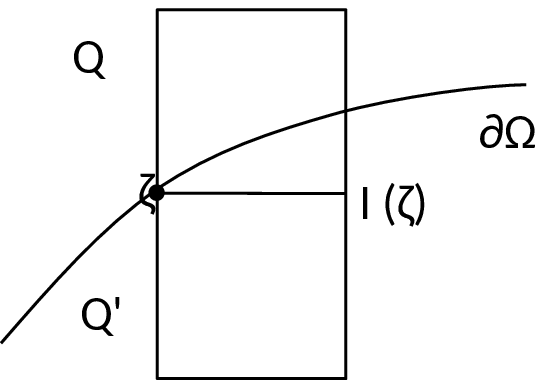}
\caption{\ } 
\label{fig_1}       
\end{figure}

So, assume that $\partial\Om$ intersects the right edge of $Q$ and $I(\za)$ is the lower edge of $Q$.
Then there exists a function $\ph\in \mathcal{C}^{1,\om}$ such that
the squares $Q(\za)$ and $Q^\prime(\za)$ with edge $I(\za)$ have the following properties:
\[
\begin{split}
&Q^\prime(\za)\subset \Om,
\\
&\Om\cap Q(\za) = \{(x,y)
= x+iy\in\Cbb: (x, y_1)\in I(\za),\ y_1\le y < \ph(x)\},
\\
&\textrm{dist}(I^\prime(\za), \partial\Om)\ge C\qle \quad \textrm{for a constant\ } C\in(0, 1),
\end{split}
\]
where $I^\prime(\za)= \{(x, y_1-\qle): (x, y_1)\in I(\za)\}$ is the lower edge of $Q^\prime(\za)$.
Using compactness of $\partial\Om$, we fix $r_0>0$ such that the above properties
hold for all $\za\in\partial\Om$ with $Q(\za)$ such that $\qle(Q(\za))=r_0$.
See Figure~\ref{fig_1}.

\subsection{Bloch spaces}\label{ss_bloch}
Given a modulus of continuity $\om$, the Bloch space $\bloch_\om(\Om)$
consists of holomorphic in $\Om$ functions $f$ such that
\[
\sup_{z\in\Om} \frac{|f^\prime(z)|\rho(z)}{\om(\rho(z))} < \infty,
\]
where $\rho(z)$ denotes the distance from $z$ to $\partial\Om$.

If the modulus of continuity $\om$ is Dini-smooth and
$\Om\subset\Cbb$ is a bounded domain with a $\ccom$-smooth boundary,
then standard arguments guarantee that
\[
\bloch_\om(\Om)\subset L^\infty(\Om) \subset L^1(\Om).
\]

The following result shows that the Beurling transform of the characteristic function
$\chi_\Om$ is in an $\om$-weighted Bloch space.

\begin{proposition}[{\cite[Theorem~2]{Va15}}]\label{p_bloch}
Let $\om$ be a regular modulus of continuity and
let $\Om\subset\Cbb$ be a bounded domain with a $\ccom$-smooth boundary.
  Then $\beu\har_{\Om}\in \bloch_{\om} (\Cbb\setminus \partial \Om)$.
\end{proposition}

In fact, one may consider Proposition~\ref{p_bloch} as a Bloch-type analog
of geometric characterizations of $\partial\Om$ by Cruz and Tolsa \cite{CT12} and Tolsa \cite{To13}.

\subsection{Integral means}
An analog of the following fact is well-known for the classical space $\bmo(\Cbb)$.

\begin{lemma}\label{l_means}
Let $\om$ be a modulus of continuity and let $g\in\bmo_\om(\Cbb)$.
 Then
 \[
 |g_Q - g_{2Q}| \lesssim \om(2\qle).
 \]
 \end{lemma}
\begin{proof}
We have
\[
\begin{split}
    |g_Q - g_{2Q}| &= \left|\frac{1}{|Q|} \int_{Q} g(z)\, dA(z) - g_{2Q}\right| \\
     &\le \frac{1}{|Q|} \int_{Q} \left| g(z)  - g_{2Q}\right|\, dA(z)  \\
     &\le \frac{1}{|Q|} \int_{2Q} \left| g(z)  - g_{2Q}\right|\, dA(z) \\
     &\le 2\om(2\qle)
\end{split}
\]
by the definitions of $g_Q$, $g_{2Q}$ and $\bmo_\om$.
\end{proof}

\subsection{An embedding lemma}
We will need the following basic observation
(see also \cite{Ja76} for similar facts related to $\bmo_\om(\Rbb)$).
\begin{lemma}\label{l_bmo_mean}
Let $\om$ be a modulus of continuity and let $f\in L^\infty(\Om)$. Assume that
\[
\sup_{Q\subset \Cbb} \frac{1}{|Q|\om(\qle)} \int_{Q\cap \Om} |f(z) - b_Q|\, dA(z)
=K< \infty
\]
for some constants $b_Q \in\Cbb$. Then $f\in \bmo_\om(\Om)$.
\end{lemma}
\begin{proof}
For any square $Q\subset \Cbb$, we have
\[
\begin{split}
   \int_{Q\cap \Om} |f(z) - f_{\qom}|
   &\, dA(z) \le  \int_{Q\cap \Om} |f(z) - b_Q|\, dA(z) + \int_{Q\cap \Om} |f_{\qom} -b_Q|\, dA(z)\\
     &\le K|Q|\om(\qle) + \int_{Q\cap \Om} \frac{1}{|Q\cap \Om|}\int_{Q\cap \Om} |f(w)-b_Q| \, dA(w)\, dA(z) \\
     &\le 2K|Q|\om(\qle).
\end{split}
\]
So, $f\in \bmo_\om(\Om)$.
\end{proof}

\begin{lemma}\label{l_emb}
Let $\om$ be a regular modulus of continuity and let $\Om\subset\Cbb$ be a bounded domain with a $\ccom$-smooth boundary.
Then
\[
\bloch_\om(\Omega) \subset \bmo_\om(\Omega).
\]
\end{lemma}
\begin{proof}
First, as indicated in Section~\ref{ss_bloch}, $\bloch_\om(\Omega)\subset L^\infty(\Om)$.
Next, to estimate the seminorm $\|Bf\|_{\bmo_\om(\Om)}$,
we have to consider the supremum over all squares $Q$, as defined by \eqref{e_camp_df}.
Let $r_0=r_0(\Om)>0$ be the constant fixed in Section~\ref{ss_r0}.
If $Q\subset\Cbb$ is a square such that $\qle(Q)\ge r_0$, then
\[
\frac{1}{|Q|}\int_{Q} |f(z)|\, dA(z) \lesssim \|f\|_{L^\infty(\Om)} \lesssim \|f\|_{\bloch_\om(\Omega)} \lesssim \om(r_0).
\]
Hence, it remains to consider the squares $Q$ with $\qle(Q)< r_0$.
If $Q\cap\Om=\varnothing$, then the required estimate is trivial.
So, let $Q\cap\Om\neq\varnothing$.
After a suitable rotation of $\Cbb$ (by $\pm\frac{\pi}{2}$ or $\pi$), we
guarantee that the lower edge of $Q$ intersects $\Om$.
Moreover, we may assume
that the lower edge of $Q$ is in $\Om$, applying a suitable rotation and
replacing $Q$ by a smaller
square such that $Q\cap\Om$ is not changed, or
shifting $Q$ down so that the intersection $Q\cap\Om$ becomes larger.
If $Q\subset\Om$, then let $(b, y_1)$ denote the upper right vertex of $Q$
else let $(b, y_1)$ denote the intersection point of $\partial\Om$
and the right edge of $Q$.

Using the properties of $Q\cap\Om$ described in the definition of $r_0$,
we conclude that
\[
Q\cap\Om = \{(x,y)=x+iy\in\Cbb: a\le x \le b,\ y_2\le y < \psi(x) \},
\]
where $y_1-\qle\le y_2\le y_1$, $b-a=\qle< r_0$ and $\psi(x)= \min\{\ph(x), y_2+\qle\}$.

The function $f$ is complex-valued; however, consideration of $\mathrm{Re\,} f$ and $\mathrm{Im\,} f$
allows to assume below in the proof that $f$ is real-valued.

By Lemma~\ref{l_bmo_mean}, it suffices to prove that
\begin{equation}\label{e_fnot_mean}
\frac{1}{|Q|}\int_{Q\cap\Om} |f - f_{y_0}| \lesssim \om(\qle),
\end{equation}
where $y_0 = y_1 - \qle$ and
\[
f_{y_0} = \frac{1}{b-a} \int_a^b f(t, y_0)\, dt.
\]

First, the function $f(t, y_0)$ is real-valued and continuous for $a\le t \le b$,
thus, $f_{y_0} = f(s, y_0)$ for certain $s\in (a, b)$.
Therefore, for any $x\in (a,b)$, we have
\begin{equation}\label{e_horizont}
|f(x, y_0) - f_{y_0}| = |f(x, y_0) - f(s, y_0)| \le \qle |\nabla f(\xi, y_0)|,\quad \textrm{where\ } \xi=\xi(x).
\end{equation}
By the definition of $r_0>0$, we have $\mathrm{dist}([a,b]\times y_0, \partial\Om)\ge C\qle$ with $0<C<1$.
Hence, applying \eqref{e_horizont}, we obtain
\begin{equation}\label{e_horiz_uni}
 |f(x, y_0) - f_{y_0}| \le \qle \frac{\om(C\qle)}{C\qle} \lesssim \om(\qle)\quad\textrm{uniformly in\ } x\in (a,b),
\end{equation}
since  $\frac{\om(t)}{t}$ is almost decreasing and $\om(t)$ is increasing.

Secondly, we have
\[
\begin{split}
    \frac{1}{|Q|}\int_{Q\cap\Om}
&|f(x, y) - f_{y_0}|\, dx dy \\
&\le
 \frac{1}{|Q|}\int_{Q\cap\Om} |f(x, y_0) - f_{y_0}|\, dx dy
 + \frac{1}{|Q|}\int_{Q\cap\Om} |f(x, y) - f(x, y_0)|\, dx dy \\
&:= E +F.
\end{split}
\]
Clearly, $E\lesssim \om(\qle)$ by \eqref{e_horiz_uni}.
Using the definition of $\bloch_\om(\Om)$ and Fubini's theorem, we obtain
the following chain of estimates:
\[
\begin{split}
F
&\lesssim\frac{1}{|Q|}\int_{Q\cap\Om} \int_{y_0}^y |\nabla f (x,t)|\, dt dy dx \\
&\lesssim\frac{1}{|Q|}\int_a^b \int_{y_1}^{\psi(x)} \int_{y_0}^y \frac{\om(\psi(x)-t)}{\psi(x)-t} \, dt dy dx\\
&=\frac{1}{|Q|}\int_a^b \int_{y_0}^{\psi(x)} \int_{t}^{\psi(x)} \frac{\om(\psi(x)-t)}{\psi(x)-t} \, dy dt dx\\
&\lesssim\frac{1}{|Q|}\int_a^b \int_{y_0}^{\psi(x)} \om(\psi(x)-t)\, dt dx\\
&\lesssim \frac{1}{\qle^2}\int_a^b 2\qle \om (2\qle)\, dx \\
&\lesssim \om(\qle),
\end{split}
\]
since $\om(t)$ is increasing and $\frac{\om(t)}{t}$ is almost decreasing.
In sum, we obtain \eqref{e_fnot_mean}. So, the proof of the lemma is finished.
\end{proof}

\subsection{An extension lemma}
The domain $\Om$ has a sufficiently smooth boundary, hence, one could try to extend the
argument from \cite{Jo80} to the $\om$-weighted case
in the setting of appropriate domains in $\Rbb^n$, $n\ge 2$.
However, 
we apply a different argument,
which uses the specifics of the complex plane.

The domain $\Om$ has a $\ccom$-smooth boundary, hence,
there exists a bilipschitz map $\tau: \Cbb\to\Cbb$
such that $\tau(\infty) = \infty$ and $\tau(\Tbb)=\partial\Om$, where $\Tbb=\{z\in\Cbb: |z|=1\}$;
see \cite[Theorem~7.10]{Po92}.
Let $\tau_{-1}$ denote the inverse of $\tau$. Then the map
\begin{equation}\label{e_ext_df}
\ext(z) = \tau\left( \frac{1}{\overline{\tau_{-1}}}\right)(z),\quad z\in\Cbb,
\end{equation}
is bilibschitz in a neighborhood of $\partial\Om$.
In the following lemma, using $\ext$,
we extend each function from $\bmo_\om(\Om)$ to an element of $\bmo_\om = \bmo_\om (\Cbb)$.

\begin{lemma}\label{l_ext}
Let $\Om\subset\Cbb$ be a bounded domain with a $\ccom$-smooth boundary
and let $\ext$ be defined by \eqref{e_ext_df}.
For $f\in \bmo_\om(\Om)$, define
\[
\begin{split}
   \widetilde{f}(z) &= f(z), \quad z\in\Om, \\
     \widetilde{f}(z) &= f (\ext(z)), \quad z\in\Cbb\setminus\Om.
\end{split}
\]
Then $\widetilde{f}\in\bmo_\om\cap L^\infty(\Cbb)$ and $\|\widetilde{f}\|_{\bmo_\om}\lesssim \|f\|_{\bmo_\om(\Om)}$.
\end{lemma}
\begin{proof}
We have $f\in \bmo_\om(\Om)\subset L^\infty(\Om)$, thus $\widetilde{f}\in L^\infty(\Cbb)$.
So, below we estimate the supremum defined by \eqref{e_bmo_df}.

The maps $\tau$ and $\tau_{-1}$ are bilipschitz on $\Cbb$,
so all distances are distorted at most in $M$ times, where $M$ is bilipschitz constant of $\tau$.
Therefore, it suffices to consider the case, where $\Om$ is the unit disk $\Dbb$ and $\ext(z)=\ext^{-1}(z) = 1/\overline{z}$.

Fix a constant $r_0\in (0, 1/4)$. Let $Q\subset \Cbb$ be a square.
If $\qle(Q)\ge r_0$, then
\[
\frac{1}{|Q|}\int_{Q} |f(z)|\, dA(z) \lesssim \|f\|_\infty \lesssim \om(r_0).
\]
So, it remains to consider the case, where $\qle(Q)\le r_0$.

\subsection*{$Q$ is in the complement of $\Dbb$}
Let $d=|w_0|>1$, where $w_0$ denotes the center of $Q$.
We have $\qle(Q)\le r_0 < 1/4$ and
\[
\frac{1}{d+2\ell} < \frac{1}{|z|} < \frac{1}{d-2\ell}, \quad z\in Q.
\]
Hence, using elementary geometric arguments, we obtain a square $\qti$
such that
$\ext(Q)\subset \qti\subset \Dbb$ and $\qle(\qti)\approx \frac{\qle}{d^2}$.

So, we have the following chain of inequalities:
\[
\begin{split}
   \frac{1}{|Q|}\int_{Q} |\fti(z) - f_{\qti}| \, dA(z)
          &= \frac{1}{|Q|}\int_{\ext(Q)} |f(z) - f_{\qti}| \frac{1}{|z|^4}\, dA(z)  \\
          &\lesssim \frac{d^4}{|Q|}\int_{\qti} |f(z) - f_{\qti}|\, dA(z)  \\
          &\lesssim \frac{1}{|\qti|}\int_{\qti} |f(z) - f_{\qti}|\, dA(z)  \\
          &\lesssim \om(\qle/d^2) \\
          &\le \om(\qle).
\end{split}
\]
\subsection*{$Q$ intersects $\Dbb$}
For $z\in Q$, we have
\[
\frac{1}{1+2\ell} < \frac{1}{|z|} < \frac{1}{1-2\ell}.
\]
Put $Q_0 = Q\cap \Dbb$ and $Q_1 = Q\cap (\Cbb\setminus \Dbb)$.
Observe that the map $\ext$ is bilipschitz in a sufficiently large neighborhood of $\partial\Dbb$. Hence, applying standard geometric
arguments, we obtain a square $\qti$ such that $\qti\supset Q_0 \cup \ext(Q_1)$
and $|\qti|\approx |Q|$.
Since $Q_0\sqcup Q_1 = Q$, we have
\[
\begin{split}
   \frac{1}{|Q|}\int_Q |\fti(z) - f_{\qti}|\, dA(z) &=\frac{1}{|Q|}\int_{Q_0} |f(z) - f_{\qti}|\, dA(z)
   + \int_{Q_1} |\fti(z) - f_{\qti}|\, dA(z)  \\
          & : = E+F.
\end{split}
\]
Clearly, $E\lesssim \om(\qle)$.
Since $\ext$ is bilipschitz on $Q_1$ with a universal distortion constant, we also obtain the following chain of estimates:
\[
\begin{split}
   F&= \frac{1}{|Q|}\int_{\ext(Q_1)} |f(z) - f_{\qti}| \frac{1}{|z|^4}\, dA(z)  \\
          &\lesssim \frac{1}{|Q|}\int_{\qti\cap\Dbb} |f(z) - f_{\qti}|\, dA(z)  \\
          &\lesssim \om(\qle).
\end{split}
\]

Finally, applying Lemma~\ref{l_bmo_mean} with $\Om=\Cbb$ and $b_Q=0$ or $b_Q= f_{\qti}$, we obtain $f\in \bmo_\om(\Cbb)$,
as required.
\end{proof}

\section{Proof of Theorem~\ref{t_main}}\label{s_main}

Let $f\in\bmo_\om(\Om)$. To estimate the seminorm $\|\beu f\|_{\bmo_\om(\Om)}$,
we have to consider the supremum over all squares $Q$, as defined by \eqref{e_camp_df}.
So, fix a square $Q\subset\Om$. Let $\fti\in\bmo_\om$ be an extension of $f$ provided by Lemma~\ref{l_ext}.
Put
\[
\begin{split}
  f_1 &= \fti_Q \chi_\Om; \\
  f_2 &= (f-\fti_Q) \chi_{2Q \cap \Om};\\
  f_3 &= (f-\fti_Q) \chi_{\Om \setminus 2Q}.
\end{split}
\]
Observe that $f= f_1 + f_2 + f_3$. So, to prove the theorem, it suffices
to show that
\[
\|\beu f_j\|_{\mo_{\qom}}:= \frac{1}{|Q|} \int_{Q\cap\Om} \left| \beu f_j - (\beu f_j)_{\qom} \right| \lesssim \om(\qle), \quad j=1,2,3.
\]

\subsection{First term}\label{ss_fterm}
By  Lemma~\ref{l_ext}, the extension $\fti$ is bounded on $\Cbb$, hence
$|\fti_Q|\lesssim C$.
Therefore,
\begin{equation}\label{e_bf1}
\|\beu f_1\|_{\mo_{\qom}} = \|\fti_Q \beu\chi_\Om\|_{\mo_{\qom}} \lesssim \|\beu\chi_\Om\|_{\mo_{\qom}}.
\end{equation}
Now, by Proposition~\ref{p_bloch} and Lemma~\ref{l_emb}, we have $\beu\chi_\Om \in \bloch_\om(\Om)\subset\bmo_\om(\Om)$,
that is, $\|B\chi_\Om\|_{\bmo_{\om}(\Om)}\lesssim C = C(\Om)$.
So, combining with \eqref{e_bf1}, we obtain
\[
\|\beu f_1\|_{\mo_{\qom}} \lesssim \|\beu\chi_\Om\|_{\mo_{\qom}} \le\om(\qle) \|\beu\chi_\Om\|_{\bmo_{\om}(\Om)} \lesssim \om(\qle),
\]
as required.

\subsection{Second term}\label{ss_sterm}
We clearly have
\[
\|\beu f_2\|_{\mo_{\qom}} \lesssim \frac{1}{|Q|} \int_{Q\cap\Om} |\beu f_2(z)|\, dA(z).
\]
Next, applying the definition of $f_2$, H\"older's inequality, boundedness of the Beurling transform on $L^2$,
a trivial integral inequality,
and the triangle inequality, we obtain the following chain of estimates:
\[
\begin{split}
   \frac{1}{|Q|}\int_{Q\cap \Om} |\beu f_2(z)|\, dA(z) &= \frac{1}{|Q|}\int_{Q\cap\Om} |\beu ((f-\fti_Q) \chi_{2Q\cap\Om}) (z)|\, dA(z) \\
     &\le \left( \frac{1}{|Q|}\int_{Q\cap\Om} |\beu ((f-\fti_Q) \chi_{2Q\cap\Om}) (z)|^2 \, dA(z) \right)^{\frac{1}{2}} \\
     &\lesssim \left( \frac{1}{|2Q|}\int_{2Q\cap\Om} |f(z) - \fti_Q|^2 \, dA(z)\right)^{\frac{1}{2}} \\
     &\le \left( \frac{1}{|2Q|} \int_{2Q} |\fti(z) - \fti_Q|^2 \, dA(z)\right)^{\frac{1}{2}} \\
     &\le |\fti_{2Q} - \fti_Q| +
      \left( \frac{1}{|2Q|} \int_{2Q} \left|\fti(z) - \fti_{2Q} \right|^2 \, dA(z)\right)^{\frac{1}{2}} \\
     &:= E + F.
\end{split}
\]
First, $E\lesssim \om(\qle)$ by Lemma~\ref{l_means}.
Secondly,
\[
F \le \om(\qle) \|\fti\|_{\bmo_{\om,2}}.
\]
As mentioned in the introduction, the seminorm in $\bmo_{\om,2}$ is equivalent to that in $\bmo_{\om,1} = \bmo_{\om}$.
Thus $F\lesssim \om (\qle)$.

\subsection{Third term}\label{ss_tterm}
Put $\fti_3 = \fti - \fti_Q$.
By the definitions of $\beu$ and $f_3$, we have
\[
\begin{split}
\pi
  &|Q| \|\beu f_3\|_{\mo_{\qom}}
  \\
  &=\int_{Q\cap\Om}
   \left| \int_{\Om\setminus 2Q} \frac{\fti_3(u)\, dA(u)}{(u-z)^2} - \frac{1}{|Q\cap\Om|} \int_{Q\cap\Om}
    \int_{\Om\setminus 2Q} \frac{\fti_3(u)\, dA(u)}{(u-w)^2}\, dA(w)
   \right|\, dA(z)
   \\
  &=\int_{Q\cap\Om} \frac{1}{|Q\cap\Om|} \int_{Q\cap\Om}\int_{\Om\setminus 2Q}
  |\fti_3(u)| \left| \frac{1}{(u-z)^2} - \frac{1}{(u-w)^2} \right|\,dA(u)\, dA(w)\, dA(z).
\end{split}
\]

Fix a point $z_0\in Q\cap\Om$. For $z, w\in Q\cap\Om$ and $u\in \Om\setminus 2Q$, we have
\[
\left| \frac{1}{(u-z)^2} - \frac{1}{(u-w)^2} \right| \lesssim \frac{|z-w|}{|u-z|^3} \lesssim \frac{\qle}{|u-z_0|^3}.
\]
Therefore,
\[
\begin{split}
  \|\beu f_3\|_{\mo_{\qom}}
  &\lesssim \frac{1}{|Q|}\int_{Q\cap\Om}\frac{1}{|Q\cap\Om|}\int_{Q\cap\Om}\int_{\Om\setminus 2Q}
   |\fti_3(u)| \frac{|z-w|}{|u-z|^3}
   \, dA(u)\, dA(w)\, dA(z)
   \\
  &\lesssim \qle \int_{\Om\setminus 2Q} \frac{|\fti(u)-\fti_Q|}{|u-z_0|^3}\, dA(u).
\end{split}
\]

Put $Q_k = \Om \cap \left(2^{k+1}Q \setminus 2^k Q \right)$, $k=1,2,\dots$. So
\[
\Om\setminus 2Q = \bigcup_{k=1}^\infty Q_k.
\]
For $u\in Q_k$, we have and $|u-z_0|\approx 2^k \qle$. Hence,
\[
\begin{split}
   \|\beu f_3\|_{\mo_{\qom}}
   &\lesssim \sum_{k=1}^\infty \frac{\qle}{(2^k \qle)^3} \int_{Q_k} |\fti(u) - \fti_Q | \, dA(u) \\
     & = \sum_{1\le k < m} + \sum_{k= m}^\infty,
\end{split}
\]
where $m$ is the smallest integer such that $2^m \qle \ge \frac{1}{2}$.
Also, the first sum is equal to zero if $m=1$.

First, $|\fti - \fti_Q |\lesssim C$ and $|Q_k|\lesssim (2^k \qle)^2$, thus
\[
\sum_{k= m}^\infty \lesssim \sum_{k= m}^\infty \frac{\qle |Q_k|}{(2^k \qle)^3}\lesssim
2^{-m}\lesssim
\om(\qle).
\]

Secondly, using a telescoping sum, definition~\eqref{e_bmo_df} and Lemma~\ref{l_means},
we obtain the following chain of estimates:
\[
\begin{split}
&\int_{Q_k} |\fti(u) - \fti_Q| \, dA(u)  \\
    &\lesssim \int_{Q_k} |\fti(u) - \fti_{2^{k}Q}|\, dA(u) + \int_{2^k Q}\sum_{n=0}^{k-1} |\fti_{2^{n+1}Q} - \fti_{2^n Q}|\, dA(u)
              \\
     &\lesssim (2^k \qle)^2 \sum_{n=1}^k \om(2^n\qle).
\end{split}
\]

Therefore,
\[
   \sum_{k=1}^{m-1} \lesssim \qle \sum_{k=1}^{m-1} \frac{(2^k \qle)^2}{(2^k \qle)^3} \sum_{n=1}^k \om(2^n\qle)
\]
So, changing the summation order, we obtain
\[
\begin{split}
 \sum_{k=1}^{m-1}   &\lesssim \qle \sum_{n=1}^{m-1} \om(2^n\qle) \sum_{k=n}^{m-1} \frac{1}{2^{k}\qle} \\
     &\lesssim \qle \sum_{n=1}^{m-1} \frac{\om(2^n\qle)}{2^n \qle} \\
     &\lesssim \qle \int_{2\qle}^{1} \frac{\om(t)}{t^2}\, dt\\
     &\lesssim \om(\qle),
\end{split}
\]
by \eqref{e_weak}.
In sum, we have
\[
\|\beu f_3\|_{\mo_{\qom}} \lesssim \om(\qle).
\]

So, the proof of the theorem is finished.

\bibliographystyle{amsplain}

\end{document}